\documentclass[12pt, a4paper]{article}
\usepackage{float}
\usepackage{natbib}
\usepackage{amsmath, amsthm, amssymb}
\usepackage{amsfonts}
\usepackage{epsf}
\usepackage{psfrag}
\usepackage{amsfonts}
\usepackage{amsmath}
\usepackage{enumitem}
\usepackage[all]{xy}
\usepackage[dvips]{graphicx}
\usepackage{epsfig}
\usepackage{epstopdf}
\usepackage{ifpdf}
\usepackage[UKenglish]{babel}
\usepackage{mathabx}
\usepackage{hyperref}
\hypersetup{
linkbordercolor={1 0 0}, 
citebordercolor={0 1 0} 
}
\newtheorem{thm}{Theorem}[section]

\newtheorem{prop}[thm]{Proposition}
\newtheorem{cor}[thm]{Corollary}
\newtheorem{defn}[thm]{Definition}
\newtheorem{examp}[thm]{Example}

\newtheorem{rmk}[thm]{Remark}
\begin{document}
\author{Faramarz Vafaee\\Michigan State University\\vafaeefa@msu.edu}
\title{Seifert surfaces distinguished by sutured Floer homology but not its Euler characteristic}
\maketitle
\begin{abstract}
In this paper we find a family of knots with trivial Alexander polynomial, and construct two non-isotopic Seifert surfaces for each member in our family. In order to distinguish the surfaces we study the sutured Floer homology invariants of the sutured manifolds obtained by cutting the knot complements along the Seifert surfaces. Our examples provide the first use of sutured Floer homology, and not merely its Euler characteristic(a classical torsion), to distinguish Seifert surfaces. Our technique uses a version of Floer homology,  called "\emph{longitude Floer homology}" in a way that enables us to bypass the computations related to the $SFH$ of the complement of a Seifert surface. 
\end{abstract}
\section{Introduction}\label{section:1}
It is known that every knot in $S^{3}$ bounds a Seifert surface. Seifert surfaces play an important role in knot theory and low dimensional topology in general. The minimum genus taken over all oriented surfaces that a knot $K$ bounds is called the genus of $K$. It is natural to wonder whether or not a given minimal genus Seifert surface for a knot is unique. To make sense of this question, we should be clear on the notion of equivalence between surfaces. We consider two surfaces $R$ and $R^{'}$ to be equivalent if there is an isotopy of $S^{3}$ taking $R$ to $R^{'}$. Fiberedness of a knot  is known as a sufficient condition for which its minimal genus Seifert surface is unique(see \cite{Burde1967}). However, there are many known examples of knots with non-isotopic Seifert surfaces.  See for instance \cite{Alford1970, Altman2011, Schubert1977, Hedden2008, Math1992, Kakimizu2005, Kobayashi1989, Lyon1974, Trotter1975}.
\\

Many classical tools in distinguishing surfaces deal with the surfaces' complements in $S^{3}$. These tools include, for instance, Seifert forms and the fundamental group of the surfaces' complements. They are quite powerful, but,
 they can potentially lead to tedious algebraic computations. There are also examples beyond the scope of classical tools. 
\\

\begin{figure}[t]
\begin{center}
\psfrag{c}{$c_{1}$}  
\psfrag{d}{$c_{2}$}    
    \includegraphics[scale=.4]{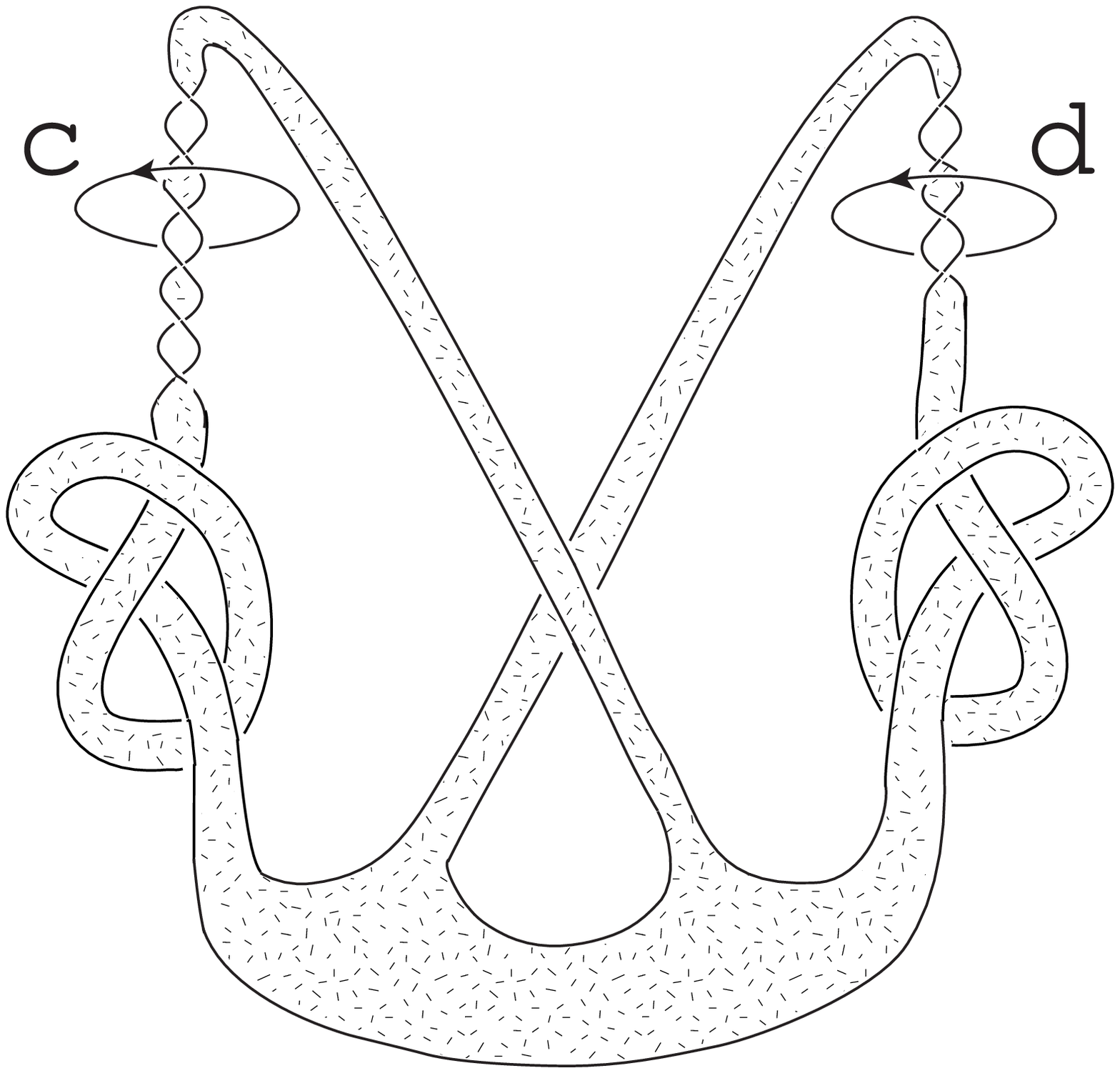}
\hfill    
\psfrag{c}{$d_{1}$}
\psfrag{d}{$d_{2}$}    
 \includegraphics[scale=.4]{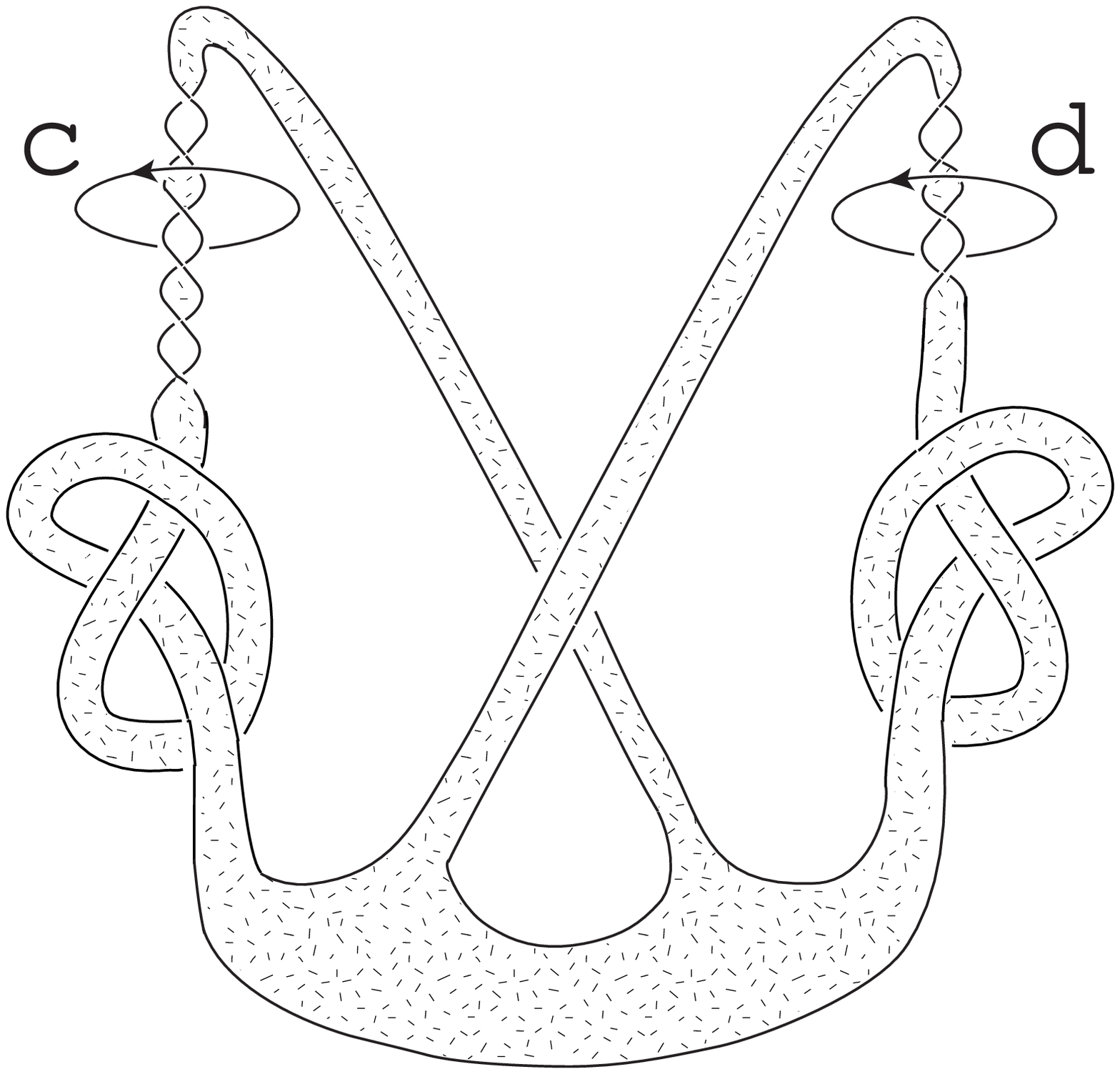}
\hfill
         \caption{The above pictures are over/under plumbings of two twisted annuli, $R$ and $R^{'}$ respectively, both are bounded by the same knot, $P(K_{1}, K_{2})$, where $K_{1}$ is the right handed trefoil and $K_{2}$ is the left handed trefoil. These lead to two distinguished Seifert surfaces $R$ and $R^{'}$, up to weak equivalence, for the knot $P(K_{1}, K_{2})$. The simple closed curves, $c_{1}$, $c_{2}$, $d_{1}$ and $d_{2}$ are basis elements for $H_{1}$ of the complement of these surfaces inside $S^{3}$.} 
\label{fig1}
\end{center}
\end{figure}

In this paper we find knots with trivial Alexander polynomial and two distinguished Seifert surfaces for each. The idea is that we plumb two untwisted annuli. Then we tie arbitrary nontrivial knots, $K_{1}$ and $K_{2}$ in each of the annuli. We produce some twists in each annulus in such a way that the framings are $l$ and 0, respectively, where $l$ is an arbitrary non-zero integer. 
We will see in section \ref{section:3} that $R$ and its dual, $R^{'}$, both are bounded by the same knot $P(K_{1}, K_{2})$. Figure \ref{fig1} shows an example. 
Our main theorem shows that the two surfaces are inequivalent, provided that one of the twisting parameters is zero and the other is non-zero. The above strategy works regardless of the knots one ties to these annuli, i.e., if you tie arbitrary knots to the annuli, the two different plumbings result in the same knot.
\\

We are now in a position to state the main theorem.
{\newtheorem*{theorema}{Theorem 1} \begin{theorema}Let $P(K_{1} , K_{2})$ be the knot obtained by plumbing two annuli with arbitrary knots $K_{1}$ and $K_{2}$ as in Figure \ref{fig1}, with framings $l$ and 0, respectively, $l \ne 0$. Changing the plumbing results in the same knot, but two inequivalent Seifert surfaces, $R$ and $R^{'}$.\end{theorema}}
Our technique begins by noting that the surfaces' complements have a particular structure called a sutured manifold(see \cite{Gabai1983}). Assigned to sutured manifolds there is an invariant called sutured Floer homology(denoted by $SFH$) introduced by  Juh\'asz in \cite{Juhasz2006}. One cannot possibly use only the rank of $SFH$ of minimal genus Seifert surfaces' complements to distinguish them, since the rank in this case depends only on the knot(see \cite[Theorem 1.5]{Juhasz2008}). Therefore, we would need to know the structure of $SFH$ as a Spin$^{c}$-graded group if we ever were able to use it to know two surfaces are not equivalent. Combining the sutured Floer homology of $(S^{3}(R), \gamma)$ with the Seifert form turns out to be a useful tool in distinguishing different Seifert surfaces(see \cite{Hedden2008}). Altman in \cite{Altman2011} gives an example of using only the sutured Floer homology polytope to distinguish two Seifert surfaces for a knot(for related definitions see Section \ref{section:2}). 
\\

The reason we are interested in the particular knots here is twofold. First, the classical methods fail in distinguishing the two Seifert surfaces. Second, the polytopes of the surfaces' complements are the same. Our theorem produces the first examples where the ($SFH$)+(Seifert form) technique is successful, but $\chi (SFH)$ alone wouldn't have sufficed. Indeed, anytime the twisting of one of the annuli is zero, we have $\chi = 0$. We refer the reader to \cite{Juhasz2010} for a detailed discussion about the identification of the Euler characteristic of the sutured Floer homology with a type of Turaev torsion polynomial.
\\

We close this section by mentioning that there are other notions of equivalence one could consider. The one we work with throughout the paper(so called \emph{weak equivalence}) is the same as regarding two surfaces $R$ and $R^{'}$ to be equivalent if there is an orientation preserving diffeomorphism between the pairs $(S^{3}, R)$ and $(S^{3}, R^{'})$(see \cite{Hatcher1983}). There is also a more restrictive notion called \emph{strong equivalence}, that considers two Seifert surfaces for a knot $K$ as equivalent if they are ambient isotopic to each other in $S^{3} \setminus n(K)$, where $n(K)$ is a neighborhood of the knot $K$ inside $S^{3}$. While it was known that the Seifert surfaces we construct via over/under plumbings for the knots  in our examples produce two distinguished Seifert surfaces for the knots up to strong equivalence(see \cite[Corollary 3.2]{Gabai1986}), what we show in the paper is stronger, that is, the two Seifert surfaces are not weakly equivalent. From now on, we will not make any further references to strong equivalence. 
\subsection*{Acknowledgements}
I would like to thank my advisor Matthew Hedden for all his invaluable support and instructive comments during the course of this work. I am also grateful to Eaman Eftekhary, Chuck Livingston, Luke Williams and David Krcatovic for helpful discussions and insights.  
\section {Background}\label{section:2}
Sutured manifolds were introduced by Gabai in \cite{Gabai1983}. Sutured Floer homology is a generalization of Ozsv\'ath-Szab\'o Floer homology to an invariant of sutured manifolds, and is defined in \cite{Juhasz2006}.
\\

In this section we begin by briefly recalling some basic notions about sutured manifolds. We then discuss the structure of $SFH$ as a group and how it behaves under decomposing a sutured manifold along an embedded surface. A key to understand it is a generalization of the Thurston norm to an invariant of relative homology classes in a sutured manifold called the Sutured-Thurston norm. We end this section by recalling a fact about the behavior of $SFH$ under a Murasugi sum of two surfaces.
\subsection{Sutured manifolds} \label{section:2.0}
The following with more details is contained in  \cite{Juhasz2006}.
{\defn \label{defn:1}  A sutured manifold $(M, \gamma)$ is a compact oriented 3-manifold with boundary, together with a set $\gamma \subset \partial M$ consisting of annuli $A(\gamma)$ and tori $T(\gamma)$. Furthermore, the interior of each component of  $A(\gamma)$ contains a suture, i.e., a homologically non-trivial simple closed curve. The union of the sutures is denoted by $s(\gamma)$.}
\\

We then take $R(\gamma) = \partial M \setminus int(\gamma)$. Define $R_{+}(\gamma)$  ($R_{-}(\gamma)$) to be those components of $R(\gamma)$ whose normal vector points out of(into) $M$. The orientation of $R(\gamma)$ must be coherent with respect to $s(\gamma)$, i.e., if $\delta$ is a component of $\partial R(\gamma)$ and is given the boundary orientation, then $\delta$ must represent the same homology class in $H_{1}(\gamma)$ as some suture. $(M, \gamma)$ is called balanced if $M$ has no closed components, $\chi(R_{-}(\gamma)) = \chi(R_{+}(\gamma))$, and the map from $\pi_{0}(A(\gamma))$ to $\pi_{0}(\partial M)$ is surjective.
\\

Let $S^{3}(R_{i}) = S^{3} \setminus int(R_{i} \times I)$ denote the complement of our Seifert surface. We equip this with a suture $\gamma = \partial R_{i} \times {\lbrace1/2\rbrace}$. 
{\defn \label{defn:2} A sutured manifold $(M, \gamma)$ is called taut if $M$ is irreducible and $R(\gamma)$ is incompressible and Thurston norm minimizing in its homology class in $H_{2}(M, \gamma)$. }
\\

In general, for a Seifert surface $R$ of a knot $K$ it follows that $S^{3}(R)$ is taut if and only if $g(R) = g(K)$. In particular, both $S^{3}(R)$ and $S^{3}(R^{'})$ in our examples are taut sutured manifolds.\\
\subsection{Sutured Floer homology polytope and decomposition of a sutured manifold along an embedded surface} \label{section:2.1}In this subsection we consider the $SFH$ of a sutured manifold as a Spin$^{c}$-graded group. We recall the definition of the sutured Floer homology polytope and then describe how the shape of the polytope changes when one decomposes a sutured manifold along an embedded surface. Throughout, we use the notation of \cite{Juhasz2010}. We do not define all the terms here and refer the reader to \cite{Juhasz2006}, \cite{Juhasz2008} and \cite{Juhasz2010} for related definitions and also more details.
\\

Let $(M, \gamma)$ be a balanced sutured manifold. Let also $v_{0}$ be a nowhere vanishing vector field pointing into $M$ along $R_{+}(\gamma)$, pointing out of $M$ along $R_{-}(\gamma)$ which restricts to $\gamma$ to be the gradient of a height function $s(\gamma) \times I \rightarrow I$. The space of such vector fields is contractible. Thus, it makes sense to fix a representative, $v_{0}$.
\\

To a sutured manifold $(M, \gamma)$, one can assign a Heegaard diagram $(\Sigma,  \mbox{\boldmath$\alpha$}, \mbox{\boldmath$\beta$})$ where $\Sigma$ is a compact oriented surface with boundary and $\mbox{\boldmath$\alpha$} = \lbrace \alpha_{1}, \alpha_{2}, ..., \alpha_{d} \rbrace$ and $\mbox{\boldmath$\beta$} = \lbrace \beta_{1}, \beta_{2}, ..., \beta_{d} \rbrace$ are two sets of pairwise disjoint simple closed curves in $int ( \Sigma ).$ Take $\mathbb{T}_{\mbox{\boldmath$\alpha$}} = \alpha_{1} \times \alpha_{2}\times ... \times \alpha_{d}$ and $\mathbb{T}_{\mbox{\boldmath$\beta$}} = \beta_{1} \times \beta_{2}\times ... \times \beta_{d}$ as subsets of the symplectic manifold Sym$^{g} ( \Sigma )$(see \cite{Juhasz2006} for details). To every point $\text{x} \in \mathbb{T}_{\mbox{\boldmath$\alpha$}} \cap \mathbb{T}_{\mbox{\boldmath$\beta$}}$, one can associate a relative Spin$^{c}$ structure, $\mathfrak{s}(\text{x}) \in \underline{Spin}^{c}(M, \gamma)$ as follows:
\\
First pick a Morse function which determines the Heegaard diagram whose gradient vector field agrees with $v_{0}$ along $\partial M$. Next, modify the vector field in a neighborhood of the flowlines specified by $\text{x}$. This produces a non-vanishing vector field $v$ that agrees with $v_{0}$ on $\partial M$. The homology class of $v$ specifies a relative Spin$^{c}$ structure which we denote by $\mathfrak{s}(\text{x})$. It turns out that $\underline{Spin}^{c}(M, \gamma)$ is an affine space over $H^{2}(M, \partial M; \mathbb{Z})$. Therefore, it makes sense to talk about the difference of two relative Spin$^{c}$ structure, $\mathfrak{s}(\text{x}) - \mathfrak{s}(\text{y}) \in H^{2}(M, \partial M; \mathbb{Z})$. We denote $PD^{-1}[\mathfrak{s}(\text{x}) - \mathfrak{s}(\text{y})]$ by $\epsilon(\text{x} , \text{y})$ which is an element of $H_{1}(M; \mathbb{Z})$.
{\defn \label{defn:3} Let $(M, \gamma)$ be a balanced sutured manifold. The support of sutured Floer homology of $(M, \gamma)$ is
\[ S(M, \gamma) = \lbrace \mathfrak{s} \in Spin^{c}(M, \gamma) : SFH(M, \gamma, \mathfrak{s}) \ne 0 \rbrace. \] }\\
$SFH(M, \gamma)$ is a finitely generated abelian group(see \cite{Juhasz2006}), thus, $S(M, \gamma)$ is finite. Moreover, if $(M, \gamma)$ is taut then $S(M, \gamma) \ne \emptyset$ by \cite[Theorem 1.4]{Juhasz2008}.
\\

It turns out that $v_{0}^{\perp}$ is a trivial vector bundle over $\partial M$, provided that $(M, \gamma)$ is balanced in each component. Let us denote the set of all trivializations of $v_{0}^{\perp}$ by $T(M, \gamma)$. For $\mathfrak{t} \in T(M, \gamma)$, let $c_{1}(\mathfrak{s}, \mathfrak{t}) \in  H^{2}(M, \partial M; \mathbb{Z})$ be the relative Euler class of the vector bundle $v^{\perp}$ with respect to the trivialization $\mathfrak{t}$, where $v$ is a nowhere zero vector field along $M$ which agrees with $v_{0}$ on $\partial M$. In other words, $c_{1}(\mathfrak{s}, \mathfrak{t})$ is the obstruction to extending $\mathfrak{t}$ from $\partial M$ to a trivialization of $v^{\perp}$ over $M$.
{\defn \label{defn:4}Fix $\mathfrak{t} \in T(M, \gamma)$. Define
\[  C(M, \gamma, \mathfrak{t}) =  \lbrace i(c_{1}(\mathfrak{s}, \mathfrak{t})) : \mathfrak{s} \in S(M, \gamma) \rbrace \subset H^{2}(M, \partial M ; \mathbb{R}), \]
where $i : H^{2}(M, \partial M; \mathbb{Z}) \rightarrow  H^{2}(M, \partial M; \mathbb{R})$ is the map induced by the natural embedding $\mathbb{Z} \hookrightarrow \mathbb{R}$. }
\\

Let $P(M, \gamma, \mathfrak{t})$ be the polytope obtained as the convex hull of $C(M, \gamma, \mathfrak{t})$ inside $H^{2}(M, \gamma; \mathbb{R})$. Thus if $(M, \gamma)$ is taut and $\alpha \in H^{2}(M, \partial M)$, then,
\[ c(\alpha, \mathfrak{t}) = min \lbrace \langle c, \alpha \rangle : c \in C(M, \gamma, \mathfrak{t}) \rbrace, \]
is a well-defined number.
{\defn \label{defn:5}For $\alpha \in H^{2}(M, \partial M)$, let
\[ H_{\alpha} = \lbrace x \in H^{2}(M, \partial M; \mathbb{R}) : \langle x, \alpha \rangle = c(\alpha, \mathfrak{t}) \rbrace. \]
 In addition, we take
\[ P_{\alpha} = H_{\alpha} \cap P(M, \gamma, \mathfrak{t}), \]
and we also use the notation,
\[ SFH_{\alpha}(M, \gamma) = \bigoplus_ {\lbrace \mathfrak{s} \in Spin^{c}(M, \gamma): i(c_{1}(\mathfrak{s}, \mathfrak{t})) \in P_{\alpha}(M, \gamma, \mathfrak{t}) \rbrace } SFH(M, \gamma, \mathfrak{s}). \]}
It turns out that this is independent of $\mathfrak{t}$.
\\

The following useful fact is contained in \cite[Proposition 4.13]{Juhasz2010}.
{\prop \label{prop:1}Let the sutured manifold $(M, \gamma)$ be taut and strongly balanced. Fix an element $\alpha \in H_{2}(M, \partial M).$ Then $P_{\alpha}(M, \gamma, \mathfrak{t})$ is a face of the polytope $P(M, \gamma, \mathfrak{t}).$ If $S$ is a nice decomposing surface that results in a taut decomposition $\xy(0,0) \xymatrix{(M, \gamma) \ar @{~>} (20,0) ^-{S}="(M, \gamma)" &\text{ \space \space}(M^{'}, \gamma^{'})} \endxy$and $[S] = \alpha,$ then,
\[ SFH(M^{'}, \gamma^{'}) \cong SFH_{\alpha}(M, \gamma).\]}
\subsection{Sutured-Thurston norm and depth of a sutured manifold} \label{section:2.2}
We now recall the definitions of different norms assigned to a sutured manifold and discuss how they are related. All is contained in \cite{Friedl2011}.\\
{\defn \label{defn:6}Let $(M,\gamma)$ be a sutured manifold. Given a properly embedded, compact connected oriented surface $S\subset M$, let,
\[ x^{s}(S) = max\lbrace 0, \frac{1}{2}|S\cap s(\gamma)| - \chi (S)\rbrace , \]
and extend this definition to disconnected surfaces by taking the sum over the components.\\
Now for $\alpha\in H_{2}(M,\partial M)$, take,
\[ x^{s}(\alpha) = min\lbrace x^{s}(S): \lbrack S, \partial S \rbrack = \alpha \rbrace .\] where the minimum is taken over all properly embedded surfaces $S \subset M$.
}
{\examp \label{examp:1}Let $K\subset S^{3}$ be a knot and let $(M,\gamma) = (S^{3} \setminus n(K), \gamma)$ where $n(K)$ is a neighborhood of the knot $K$ and $\gamma$ consists of two meridianal sutures. If $\alpha \in H_{2}(M, \partial M)$ is a generator, then $x^{s}(\alpha) = 2g(K)$. Note that this differs from the usual Thurston norm $x$ of $M$, which satisfies $x(\alpha) = 2g(K) - 1$ for a nontrivial knot.
}
{\defn \label{defn:7}Let $S(M, \gamma)$ be the support of $SFH(M, \gamma)$. If $ \alpha \in H_{2}(M, \partial M ; \mathbb{R})$, we define, 
\[z(\alpha) = max\lbrace \langle \mathfrak{s} - \mathfrak{t}, \alpha\rangle : \mathfrak{s}, \mathfrak{t} \in S(M, \gamma)\rbrace.\]
}
The following key proposition helps us to bypass the computations of the $SFH$ of the complements of our annuli.
{\prop\cite[Proposition 7.7]{Friedl2011} \label{prop:2}Let $(M,\gamma)$ be an irreducible balanced sutured manifold such that all boundary components of $M$ are tori. Then $z = x^{s}$.
}
\\

We do not aim to give a precise proof of the proposition here; however, the proof uses a key fact which is helpful during the course of computations we will do. The detailed proof is contained in \cite[ Proposition 7.7]{Friedl2011}.
\begin{proof} [Sketch of the proof.] It turns out that we can assume each of the components of $\partial M$  consists of two sutures. For $(M, \gamma)$, there is a link, $L$, in a 3-manifold $Y$, where $Y$ is obtained by Dehn filling $\partial M$ such that the $\mu _{i}$'s are meridians of the tori. For each $\mathfrak{s} \in Spin^{c}(M,\gamma)$ we obtain a relative first Chern class $c_{1}(\mathfrak{s})\in H^{2}(M, \partial M)$ such that the set $\lbrace c_{1}(\mathfrak{s}): \mathfrak{s}\in S(M, \gamma)\rbrace$ is symmetric about the origin. Then, due to \cite[Remark 8.5]{Juhasz2008}, for every $h \in H_{2}(M, \partial M)$,
\begin{eqnarray}
\label{eqn:1}
max\lbrace \langle c_{1}(\mathfrak{s}), h\rangle : \mathfrak{s} \in S(M,\gamma)\rbrace = x(h) + \sum_{i = 1}^{l}\mid \langle h, \mu _{i}\rangle  \mid.
\end{eqnarray}
Since the image of $S(M, \gamma)$ is symmetric and Spin$^{c}(M,\gamma)$ is an affine space over $H^{2}(M, \partial M)$, this is equivalent to saying 
\begin{eqnarray}
\label{eqn:2}
max\lbrace \langle \mathfrak{s} - \mathfrak{t}, h \rangle : \mathfrak{s}, \mathfrak{t} \in S(M,\gamma) \rbrace = x(h) + \sum_{i=1}^{l}\mid \langle h, \mu _{i}\rangle  \mid.
\end{eqnarray}
Note that the left sides of (\ref{eqn:1}) and (\ref{eqn:2}) must be the same which completes the proof.
\end{proof} 
\subsection{$SFH$ of the Murasugi sum of two manifolds} \label{section:2.3} In \cite{Juhasz2006}, Juh\'asz found a formula that governs the behavior of $SFH$ under a plumbing of two annuli. We recall \cite[Remark 10.8]{Juhasz2006} as a proposition here.
{\prop \label{prop:4}If a surface $R$ is a Murasugi sum of two subsurfaces $R_{1}$ and $R_{2}$, then over any field $\mathbb{F}$, we have,
\[SFH \left( S^{3}(R); \mathbb{F} \right)  \cong SFH\left( S^{3}(R_{1}); \mathbb{F} \right) \otimes  SFH\left( S^{3}(R_{2}); \mathbb{F} \right).\]  }
\\
The above formula is an isomorphism of Spin$^{c}$-graded groups. Note also that a plumbing is a special case of a Murasugi sum(see also \cite{Ni2006}). As a matter of fact, what we have in Figure \ref{fig1} is a Murasugi sum of two annuli along a 4-gon, so this proposition enables us to compute the $SFH$ of the complement of each of those annuli and then simply take their tensor product. 
\section{Using $SFH(Y(R))$ to distinguish Seifert surfaces} \label{section:3}
This section will be devoted to using the sutured Floer homology invariants in order to distinguish inequivalent Seifert surfaces of knots in the form of Figure \ref{fig1}. We keep using the same notation as in Section \ref{section:1}.  First, we show that, different plumbings of the knotted annuli have the same boundaries. Then, we will prove, the over/under plumbings lead to two distinguished Seifert surfaces for their common boundary.
\\

\subsection{$R$ and $R^{'}$ bound the same knot.}\label{section:3.0}
\begin{figure}[t]
\psfrag{K}{\scriptsize$K_{1}$}
     \includegraphics[scale=.2]{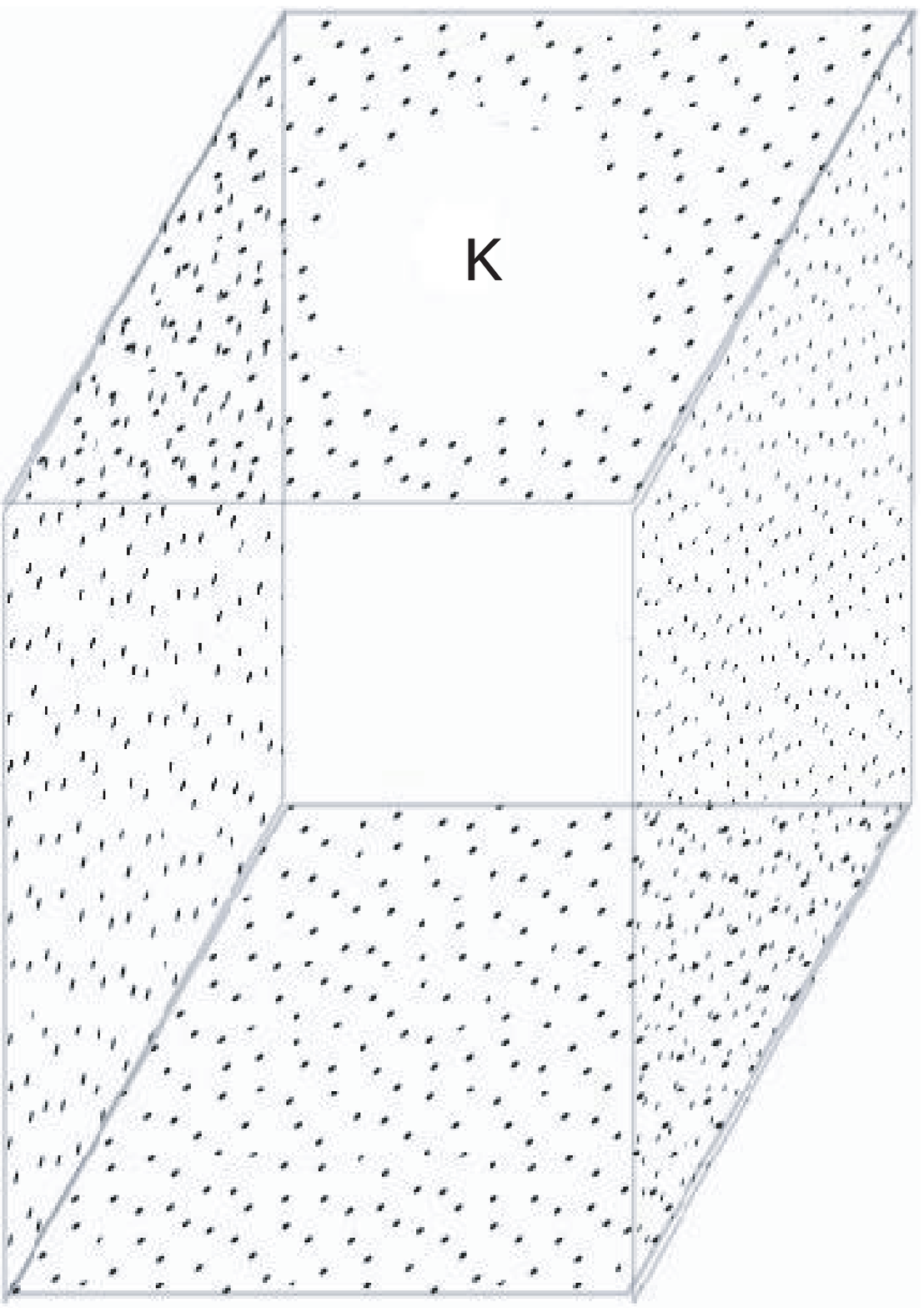}
\hfill  
\psfrag{J}{\scriptsize$K_{2}$}
     \includegraphics[scale=.2]{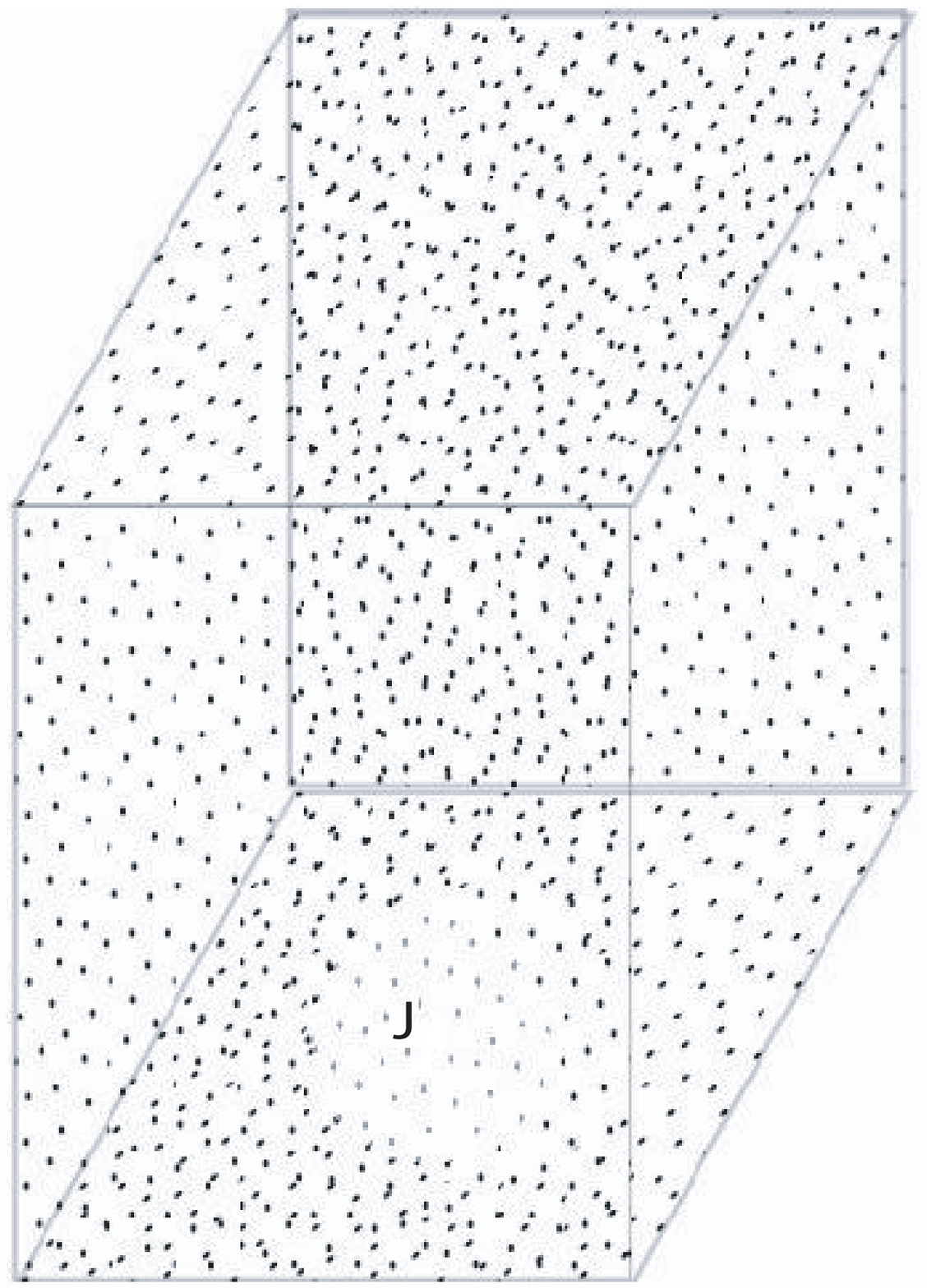}
\hfill  
\psfrag{K}{\scriptsize$K_{1}$}
\psfrag{J}{\scriptsize$K_{2}$}
      \includegraphics[scale=.23]{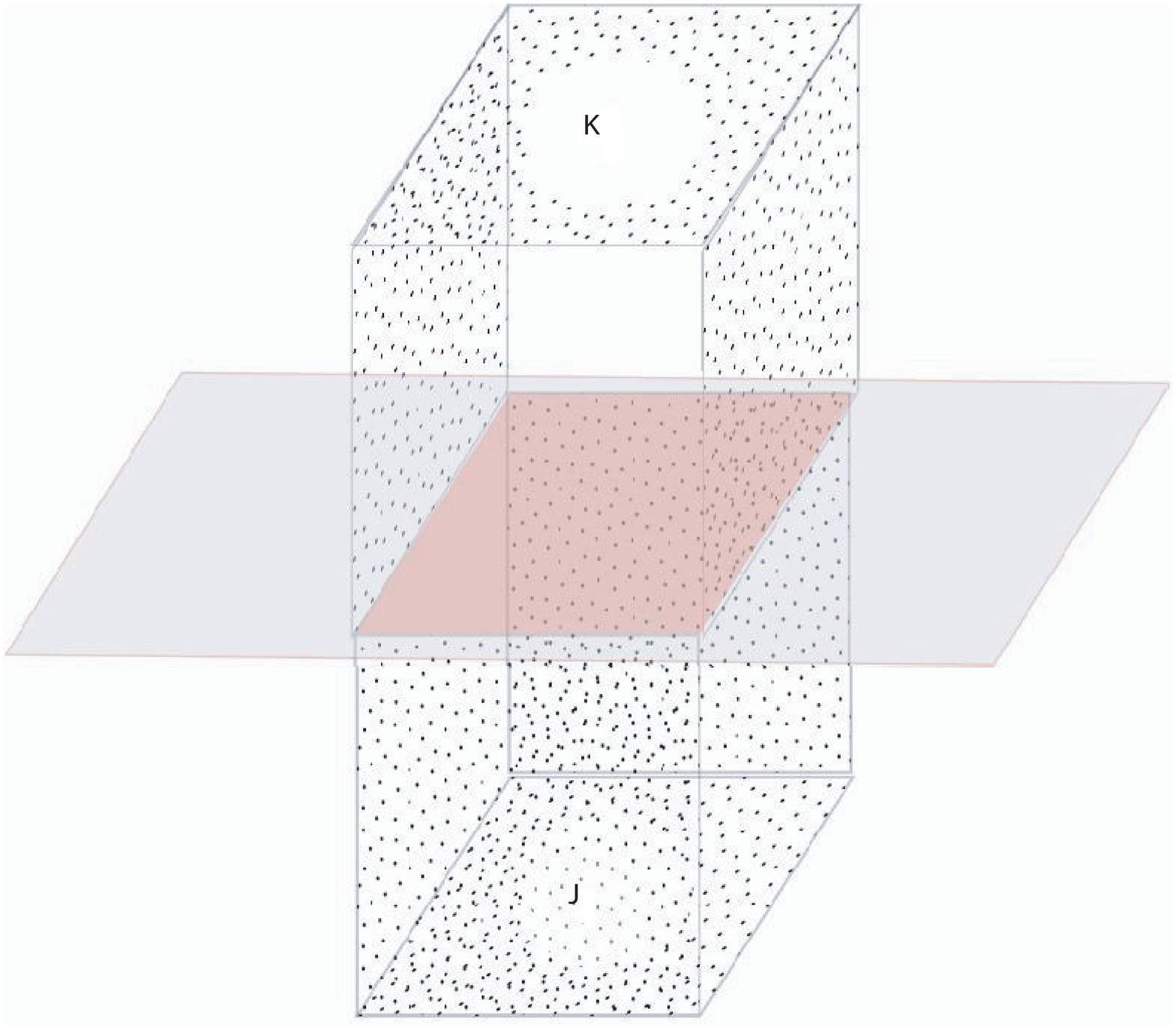}
\hfill
                \caption{The first two pictures are showing two knotted annuli. The third picture shows they are plumbed along a red 4-gon where the 4-gon is on the sphere, one annulus is inside the sphere and the other annulus is outside of it. The light blue plane is part of the sphere.} 
\label{fig2}   
\end{figure}

The goal of this subsection is to prove that the  knotted annuli in the form of Figure \ref{fig1}, where the right handed trefoil and the left handed trefoil are replaced by arbitrary knots, $K_{1}$ and $K_{2}$, both are bounded by the same knot, $P(K_{1}, K_{2})$. We would like to mention here that, it is possible to do isotopies in Figure \ref{fig1} to go from one presentation of $P(K_{1}, K_{2})$ to the other, however, we take a different route. Let us recall the definition of a plumbing here(see \cite{Gabai1983a}, \cite{Gabai1986}, \cite{Kakimizu2005} and \cite{Kobayashi1989} for more details). Given compact oriented surfaces $S_{1}, S_{2} \subset S^{3}$; if there are 3-balls $V_{1}, V_{2} \subset S^{3}$ satisfying the following properties:
\[ \begin{array}{ccc}
V_{1} \cup V_{2} = S^{3}, & V_{1} \cap V_{2} = \partial V_{1} = \partial V_{2} = S^{2}, & S_{i} \subset V_{i} \text{   } (i = 1, 2), \\
R = S_{1} \cup S_{2} \text{                                        }\text{ and} &  D =S_{1} \cap S_{2} \text{   is a 4-gon}, &
\end{array} \]  
then $R$ is called a plumbing of $S_{1}$ and $S_{2}$. In our examples $S_{1}$, $S_{2}$ and their plumbing are shown in Figure \ref{fig2}. Put $P(K_{1}, K_{2}) = \partial R$. Note that $R^{'} = (R - D) \cup D^{'}$ is an oriented surface with $\partial R^{'} = P(K_{1}, K_{2})$ where $D^{'} = S^{2} - int(D).$ We will say that $R^{'}$ is a \emph{dual} of $R$. Notice that $R^{'}$ is also a plumbing of $S_{1}$ and $S_{2}$ where $S^{'}_{i} = (S_{i} - D) \cup D^{'}$ $(i = 1, 2)$.\\
\clearpage
\begin{figure}[t]
\psfrag{K}{\scriptsize$K_{1}$}
\psfrag{J}{\scriptsize$K_{2}$}      
\includegraphics[scale=.22]{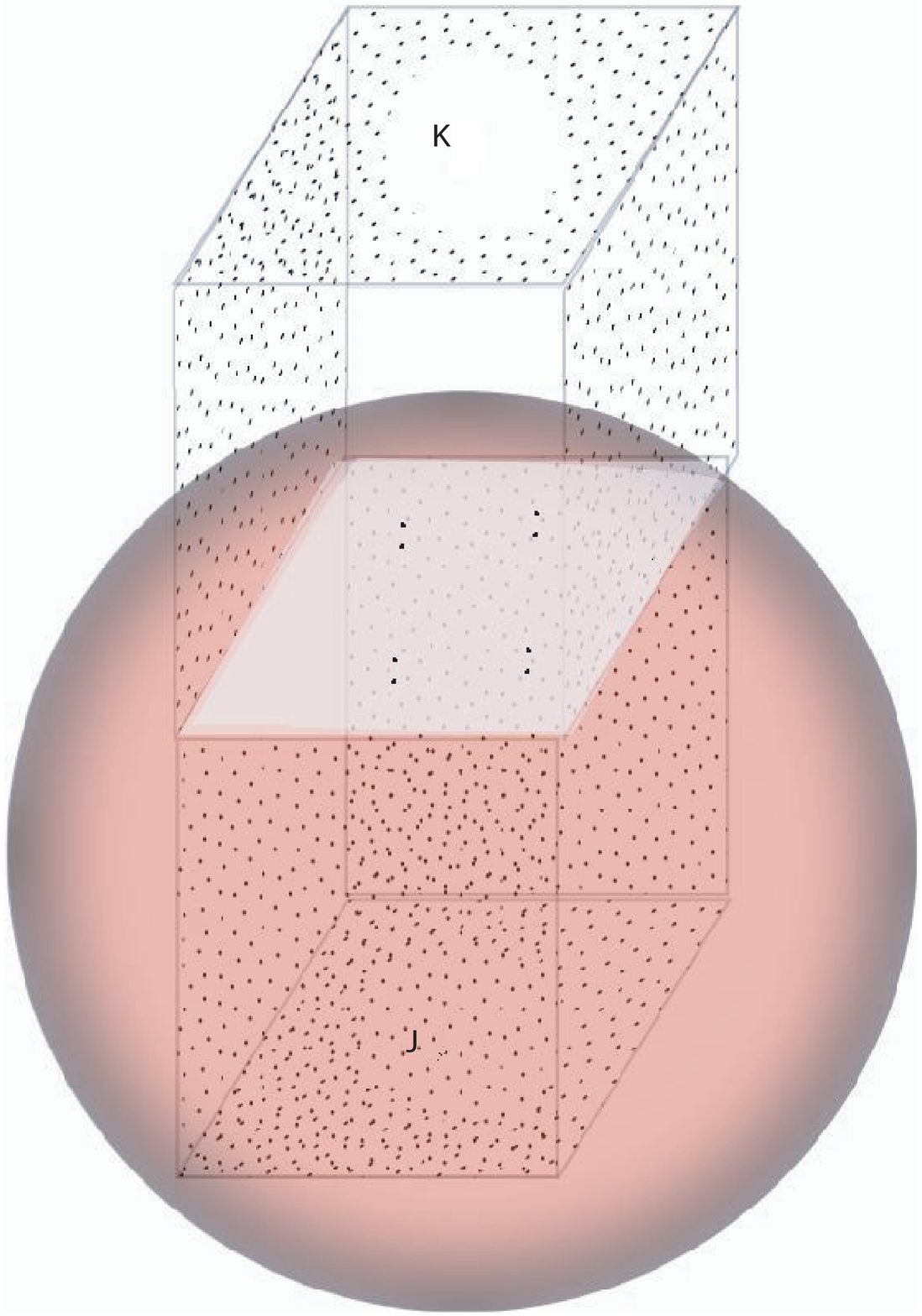}
\hfill
     \includegraphics[scale=.25]{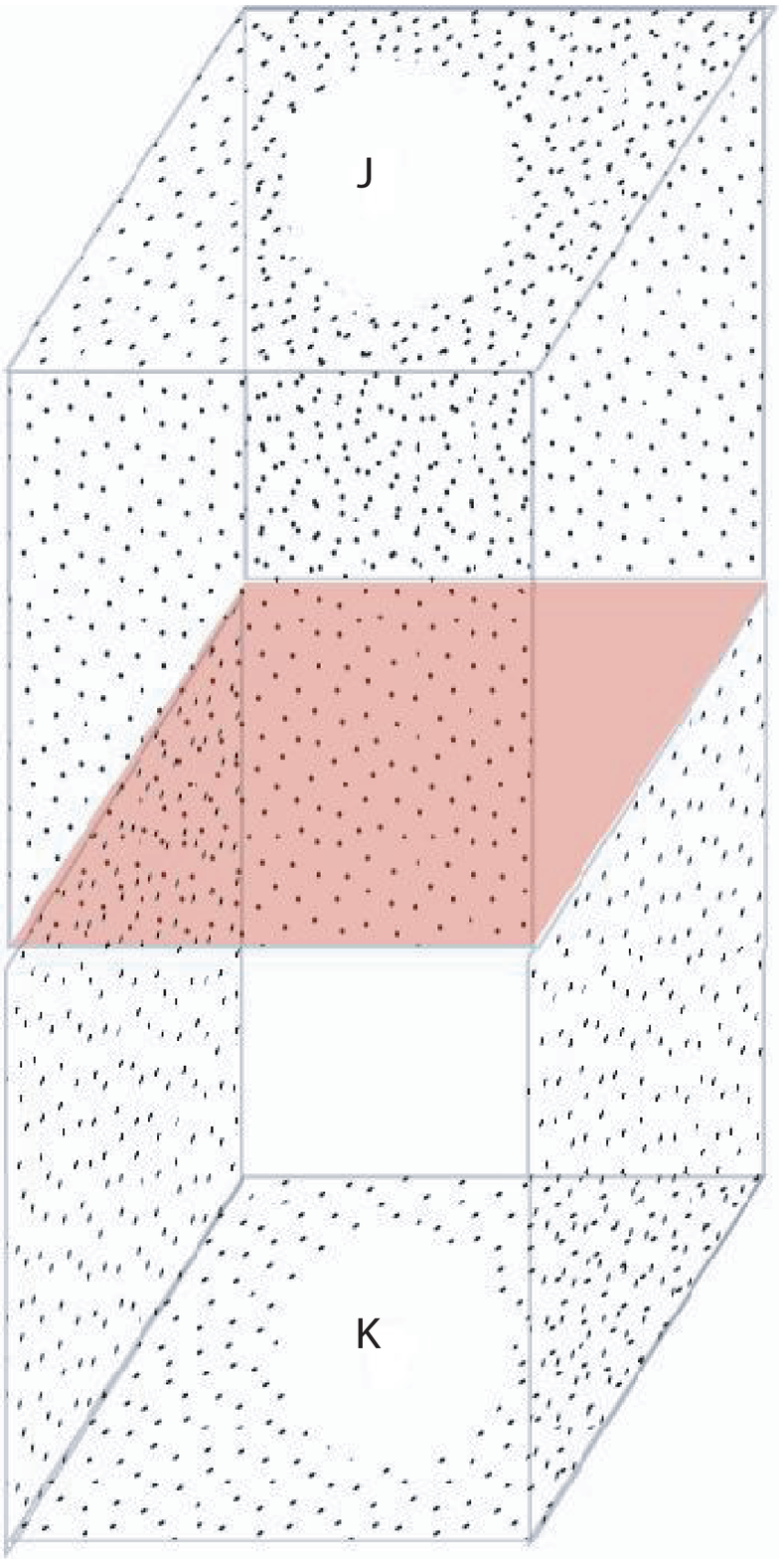}
\hfill
   \caption{ The plumbing region, $D^{'}$, is represented by the whole sphere where the 4-gon, $D$ is removed from it. To have $D^{'}$ in between the two annuli, the top annulus will get moved down and the bottom annulus will go to the top.} 
\label{fig10}
\end{figure}
\clearpage
We still need to argue that the two different plumbings we talked about above are the same as the plumbings in Figure \ref{fig1}. Notice that if we pull down the $K_{2}$-knotted annulus in the left picture of Figure \ref{fig1}, we obtain the representation of a plumbing in the sense of Figure \ref{fig2}. Then, with the above notation, if we change the plumbing region, $D$ to $D^{'}$, we will have a hole in the sphere as in Figure \ref{fig10}. Now, in order to have the red plumbing region in Figure \ref{fig10}, between the knotted annuli, we need to push the $K_{1}$-knotted annulus down, and the $K_{2}$-knotted annulus will get moved up. But, the latter is the same as the right side of Figure \ref{fig1}. 
{\rmk The above argument essentially shows that if we push the Seifert surfaces $R$ and $R^{'}$ inside the four-ball and keep the knot, $P(K_{1}, K_{2})$, in $S^{3}$, they are isotopic. Pushing the plumbing disk, $D$, of the surface $R$ inside the four-ball, enables us to go from the right picture in Figure \ref{fig2} to the right picture in Figure \ref{fig10} and simultaneously avoid all intersections that could possibly occur. Therefore, we obtain an isotopy between $R$ and $R^{'}$ inside the four-ball.}
\subsection{Classical Methods} \label{section:3.1}
\begin{figure}[t]
   \begin{center}  
\psfrag{K}{\scriptsize$K_{1}$}
\psfrag{J}{\scriptsize$K_{2}$}    
      \includegraphics[scale=.28]{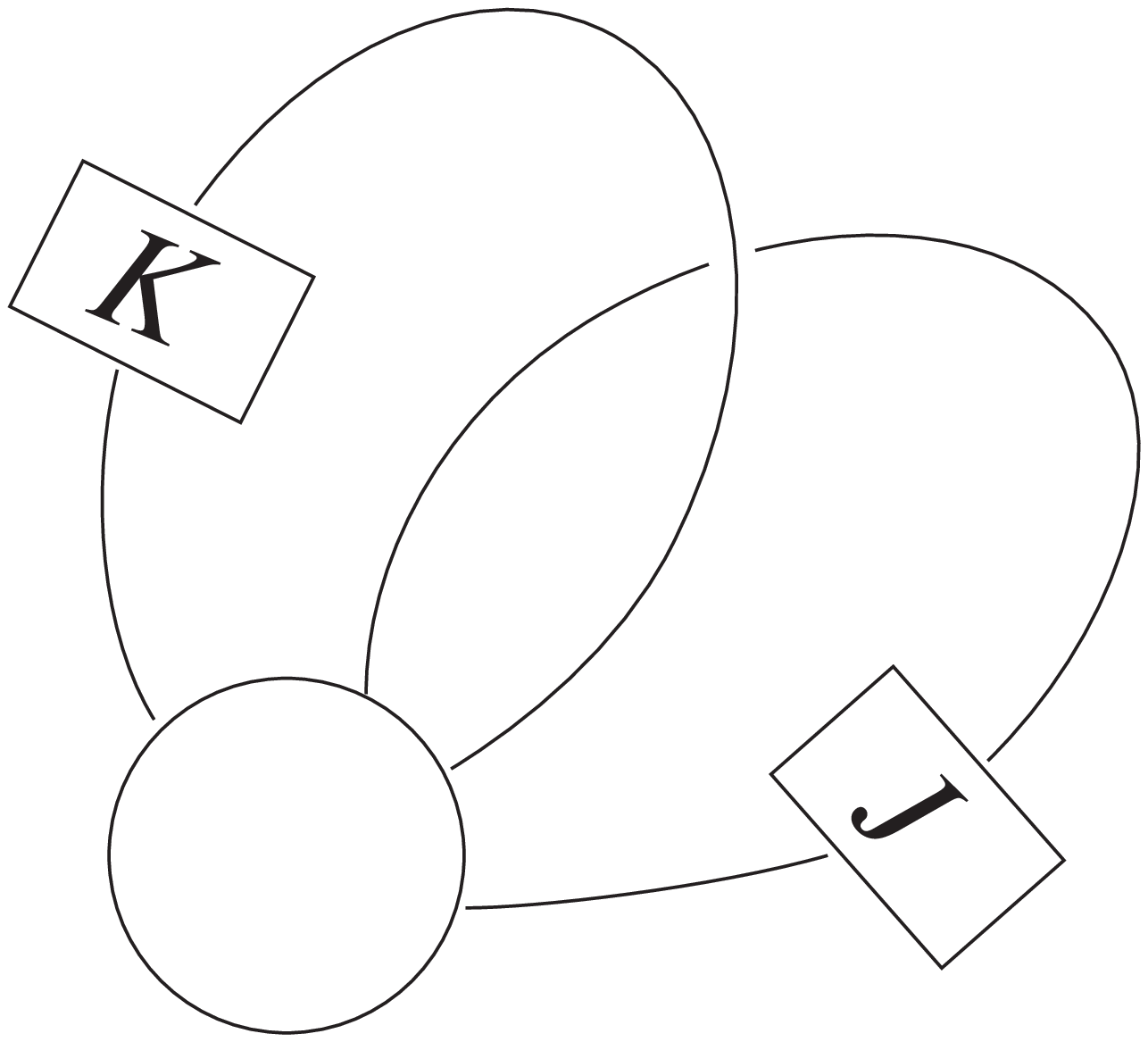}
\hfill
\psfrag{K}{\scriptsize$K_{1}$}
\psfrag{J}{\scriptsize$K_{2}$}
       \includegraphics[scale=.28]{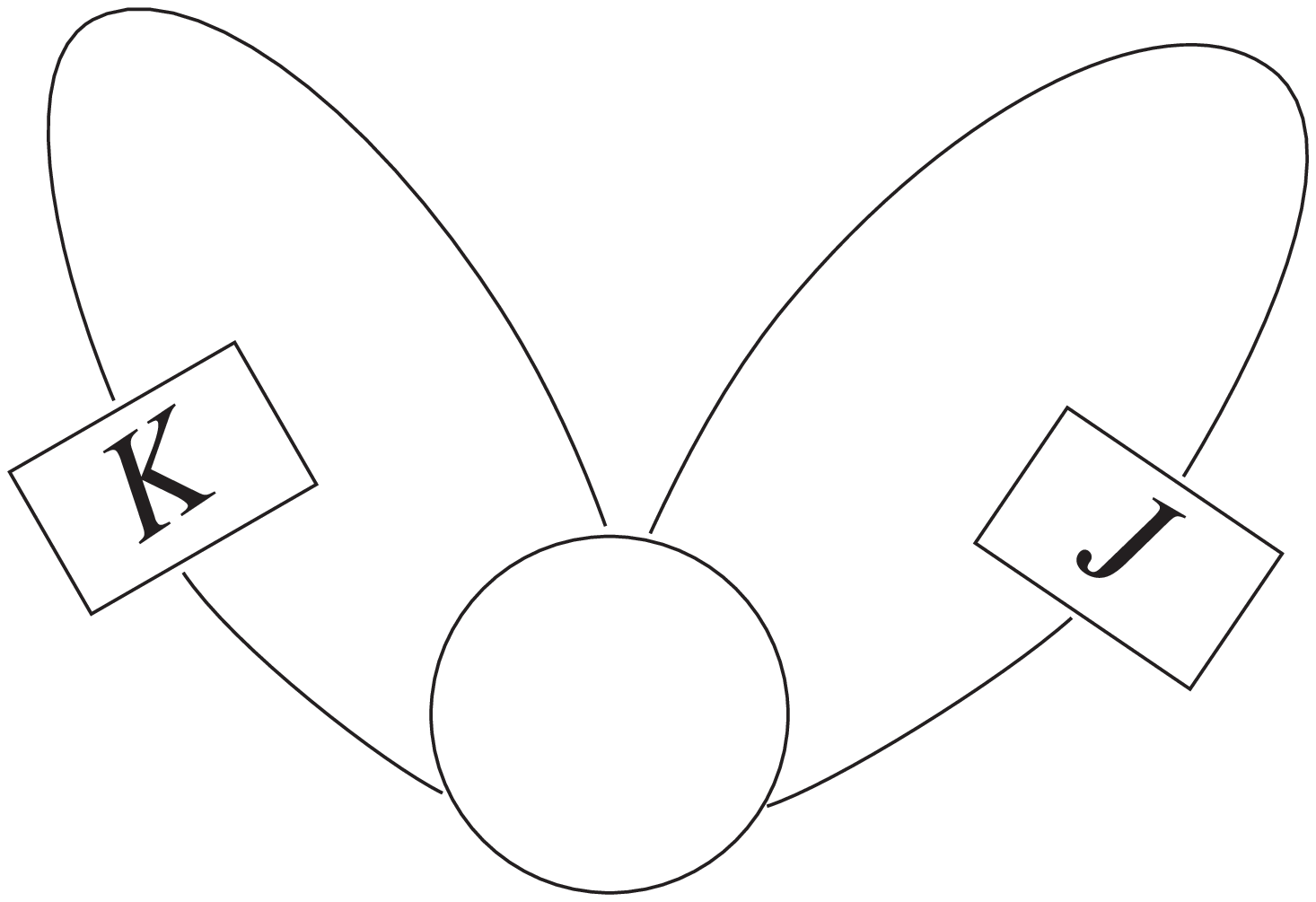}
\hfill
\psfrag{K}{\scriptsize$K_{1}$}
\psfrag{J}{\scriptsize$K_{2}$}
       \includegraphics[scale=.28]{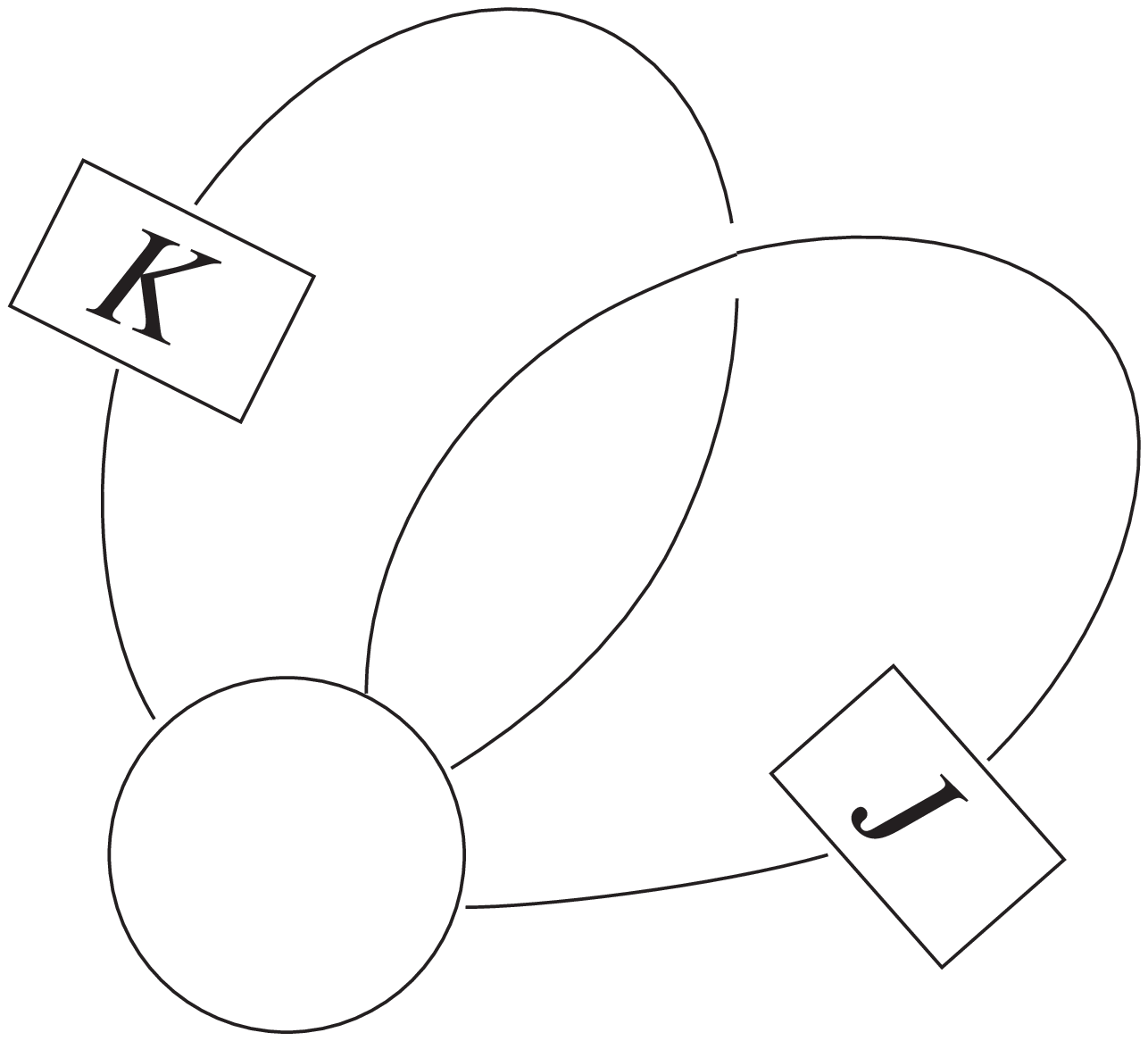}
\hfill
             \end{center}
 \caption{We encode the surfaces in this picture with planar diagrams where the circles and the bands are representing the zero cells and the one cells, respectively.} 
  \label{fig3}
\end{figure}

In this subsection we observe that classical methods fail in distinguishing the two Seifert surfaces, at least for a subfamily of our examples. Let us start by looking at the Seifert forms. Pick that subset of our examples where the non-zero twisting parameter is 1. Easy computation based on the basis elements represented in Figure \ref{fig1}, shows that the Seifert forms are given by $V_{R} = \left( \begin{array}{cc}
1 & 1 \\
0 & 0  \end{array} \right)$  
and $ V_{R^{'}} = \left( \begin{array}{cc}
1 & 0 \\
-1 & 0  \end{array} \right)$
 for $R$ and $R^{'}$, respectively. One can check that if $ W = \left( \begin{array}{cc}
 1 & -1 \\
0 & 1 \end{array} \right)$, then we have  $V_{R^{'}}$ = $W^{T}$$V_{R}$W. Since, $W \in SL_{2}(\mathbb{Z})$, $V_{R}$ and $V_{R^{'}}$ are congruent. Thus,  Seifert forms are incapable in distinguishing these particular surfaces. Another effective way to distinguish between surfaces is by looking at the homeomorphism type of their complements. Starting from the left picture in Figure \ref{fig1}, we take a regular neighborhood of the surface, $R$, inside $S^3$, we obtain a handlebody with one zero handle and two one handles. The left picture in Figure \ref{fig3} illustrates a schematic diagram for our discussion when we retract the handles to one zero cell and two one cells. Notice that we can do handle slides inside $S^3$ in the sense that in the schematic diagram, first, we obtain the middle picture in Figure \ref{fig3} and then, the right picture, that is a planar diagram for the surface, $R^{'}$. Hence, $S^{3}\setminus R$ and  $S^{3}\setminus R^{'}$ are homeomorphic. Thus any algebro-topological invariant derived from the homeomorphism type(e.g. $\pi_{1}$) will fail to distinguish the surfaces $R$ and $R^{'}$.
\subsection{Structure of $SFH(S^{3}(R))$}\label{section:3.3}As we discussed in section \ref{section:1}, we plumb two knotted annuli. Figure \ref{fig1} shows an example when the knots, $K_{1}$ and $K_{2}$, are the left handed and right handed trefoils, respectively.
The point is that the complement of each of these annuli in $S^{3}$ deformation retracts to the knot complement, regardless of what the framings of those knots are. Thus, we have $H_{1}(M
_{i}) \cong \mathbb{Z}$ where $M_{i} = S^{3}(A(K_{i}))$ and $A(K_{i})$'s are the annuli as in Figure \ref{fig4}($i = 1, 2$). Notice also that we are in a position to use Proposition \ref{prop:2}. Therefore we have $x^{s}(A(K_{i})) = z(A(K_{i})) $. We recall that 
\[ x^{s}(A(K_{i})) =  max\lbrace 0, \frac{1}{2}|A(K_{i})\cap s(\gamma)| - \chi (A(K_{i}))\rbrace \]
and that $\ x^{s}(A(K_{i})) = 2g-1$, similar to the computation we did in Example \ref{examp:1}; since $g \ge 1$, the breadth is at least 1. Thus we have the following polytopes,
\[ \setlength{\unitlength}{1mm}
    \begin{picture}(65,15)(0,-10)
      \put(0,0){\line(400,0){65}}
      \multiput(10,0)(15,0){4}{\line(0,-1){1}}
      \put(10,-4){$G_{1}$}
      \put(25,-4){$G_{2}$}
      \put(40,-4){...}
      \put(55,-4){$G_{n}$}
      \put(17, 1.4){$c_{1}$}
      \put(32, 1.4){$c_{1}$}
      \put(47, 1.4){$c_{1}$}
\end{picture} \] 
\[ \setlength{\unitlength}{1mm}
    \begin{picture}(65,15)(0,-10)
      \put(0,0){\line(400,0){65}}
      \multiput(10,0)(15,0){4}{\line(0,-1){1}}
      \put(10,-4){$H_{1}$}
      \put(25,-4){$H_{2}$}
      \put(40,-4){...}
      \put(55,-4){$H_{m}$}
      \put(17, 1.4){$c_{2}$} 
      \put(32, 1.4){$c_{2}$}
      \put(47, 1.4){$c_{2}$}
\end{picture} \]
for $S^{3}(A(K_{1}))$ and $S^{3}(A(K_{2}))$, respectively, where $G_{1}$, $G_{n}$, $H_{1}$ and $H_{m}$ are all non-zero and also we have that $m, n \ge 2$.
We now plumb the annuli, so,
\\
\\
\\
\\
\\
\\
\\
\\
\\
\\
\\
 due to Proposition \ref{prop:4} we get a tensor product formula. Thus, one obtains the polytope shown above for the manifold $S^{3}(R)$
where $c_{1}$ and  $c_{2}$ are the standard basis elements for $H_{1}(S^{3}(R); \mathbb{Z})$, specified in Figure \ref{fig1}. In addition, for $S^{3}(R^{'})$, we follow the exact same process that obviously results in the same polytope, except $c_{1}$ and $c_{2}$ are replaced by $d_{1}$ and $d_{2}$. In summary, we find that the polytopes of $SFH(S^3(R), \gamma)$ and $SFH(S^3(R^{'}), \gamma^{'})$ are rectangular and at least four of the vertices in each rectangle have non-zero groups sitting on them.   

\subsection{Calculation}\label{section:3.4} In this subsection we present a way to calculate the groups $G_{1} \otimes H_{1}$, $G_{1} \otimes H_{m}$, $G_{n} \otimes H_{1}$ and $G_{n} \otimes H_{m}$ sitting on four of the vertices of the above polytope. We explain the computations needed to obtain $G_{1}$ and $G_{n}$. Then, $H_{1}$ and $H_{m}$ could be calculated in a quite similar way. Throughout, we fix the knot $K_{1}$ and for the sake of simplicity of notation, we set $K = K_{1}$. \\

\emph{A priori} we deal with the complement of a knotted annulus, with two oriented sutures on the two edges of it(see Figure \ref{fig4}).
Let us denote this manifold by $M$. Then, $M$ is homeomorphic to the knot complement. As a sutured manifold, however, it is different; the sutures are not meridianal. They are rather like oriented longitudes.
\begin{figure}[!t]
\setlength{\unitlength}{0.8cm}
\begin{center}
\begin{picture}(8,0)
\put(0,0){\line(1,0){1.5}}
\put(3.5,0){\line(1,0){1.5}}
\put(6,0){\line(1,0){1.5}}
\put(-2,0){$G_{1} \otimes H_{1}$}
\put(1.6,0){$G_{2} \otimes H_{1}$}
\put(5.2,0){$...$}
\put(7.7,0){$G_{n} \otimes H_{1}$}
\put(.8,0.2){$c_{1}$}
\put(4.2,0.2){$c_{1}$}
\put(6.7,.2){$c_{1}$}
\put(-1,-.2){\line(0,-1){1.5}}
\put(-1,-2.5){\line(0,-1){1.5}}
\put(-1,-4.8){\line(0,-1){1.5}}
\put(-2,-2.2){$G_{1} \otimes H_{2}$}
\put(-1.1,-4.3){$.$}
\put(-1.1,-4.4){$.$}
\put(-1.1,-4.5){$.$}
\put(-2,-6.8){$G_{1} \otimes H_{m}$}
\put(-1.5,-0.9){$c_{2}$}
\put(-1.5,-3.2){$c_{2}$}
\put(-1.5,-5.5){$c_{2}$}
\put(0,-2.2){\line(1,0){1.5}}
\put(3.5,-2.2){\line(1,0){1.5}}
\put(6,-2.2){\line(1,0){1.5}}
\put(1.6,-2.2){$G_{2} \otimes H_{2}$}
\put(5.2,-2.2){$...$}
\put(7.7,-2.2){$G_{n} \otimes H_{2}$}
\put(.8,-2){$c_{1}$}
\put(4.2,-2){$c_{1}$}
\put(6.7,-2){$c_{1}$}
\put(2.5,-.2){\line(0,-1){1.5}}
\put(2.5,-2.5){\line(0,-1){1.5}}
\put(2.5,-4.8){\line(0,-1){1.5}}
\put(2.4,-4.3){$.$}
\put(2.4,-4.4){$.$}
\put(2.4,-4.5){$.$}
\put(1.55,-6.8){$G_{2} \otimes H_{m}$}
\put(2,-0.9){$c_{2}$}
\put(2,-3.2){$c_{2}$}
\put(2,-5.5){$c_{2}$}
\put(8.6,-.2){\line(0,-1){1.5}}
\put(8.6,-2.5){\line(0,-1){1.5}}
\put(8.6,-4.8){\line(0,-1){1.5}}
\put(8.5,-4.3){$.$}
\put(8.5,-4.4){$.$}
\put(8.5,-4.5){$.$}
\put(7.8,-6.8){$G_{n} \otimes H_{m}$}
\put(8.1,-0.9){$c_{2}$}
\put(8.1,-3.2){$c_{2}$}
\put(8.1,-5.5){$c_{2}$}
\put(0.1,-6.8){\line(1,0){1.4}}
\put(3.6,-6.8){\line(1,0){1.4}}
\put(6,-6.8){\line(1,0){1.5}}
\put(5.2,-6.8){$...$}
\put(.7,-6.6){$c_{1}$}
\put(4.1,-6.6){$c_{1}$}
\put(6.7,-6.6){$c_{1}$}
\end{picture}
\label{fig8}
\end{center}
\end{figure}
{\rmk \label{rmk:1} $SFH(M, \gamma)$ is isomorphic to an invariant of knots called "Longitude Floer Homology" of either the zero or $l$ framed knot. As in the construction of knot Floer homology(\cite{Ozsvath2004, Rasmussen2003}) we first find a Heegaard diagram for the knot complement in $S^{3}$. While in there we add the meridian of the knot to the set of $\beta$ curves to obtain a Heegaard diagram for $S^{3}$; in the longitude Floer homology case, we add a longitude that results in a Heegaard diagram for $S^{3}_{0}(K)$.
This subject has been first studied by Eftekhary in \cite{Eftekhary2005} for framing zero and later, has been developed by Hedden for arbitrary surgery coefficients in \cite{Hedden2007}. }
\\

Let $S(K)$ be a Seifert surface for the knot $K$. The key to obtaining $G_{1}$ and $G_{n}$ is to decompose $(M, \gamma)$ along $S(K)$. The following proposition can be enlightening.
\begin{figure}[t]
    \begin{center}  
\psfrag{K}{$K_{1}$}
      \includegraphics[scale=.3]{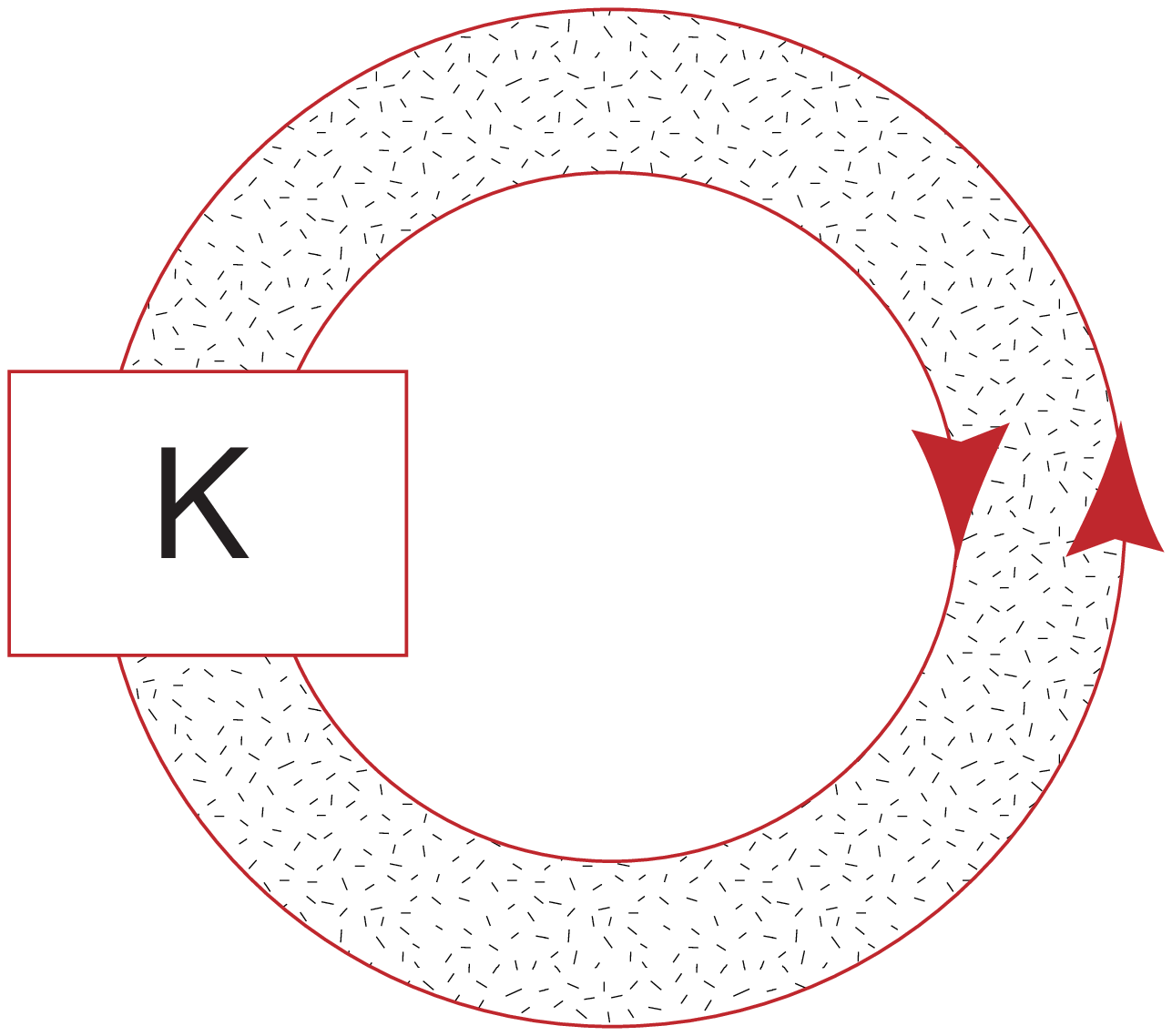}
\hfill
\psfrag{J}{$K_{2}$}
       \includegraphics[scale=.3]{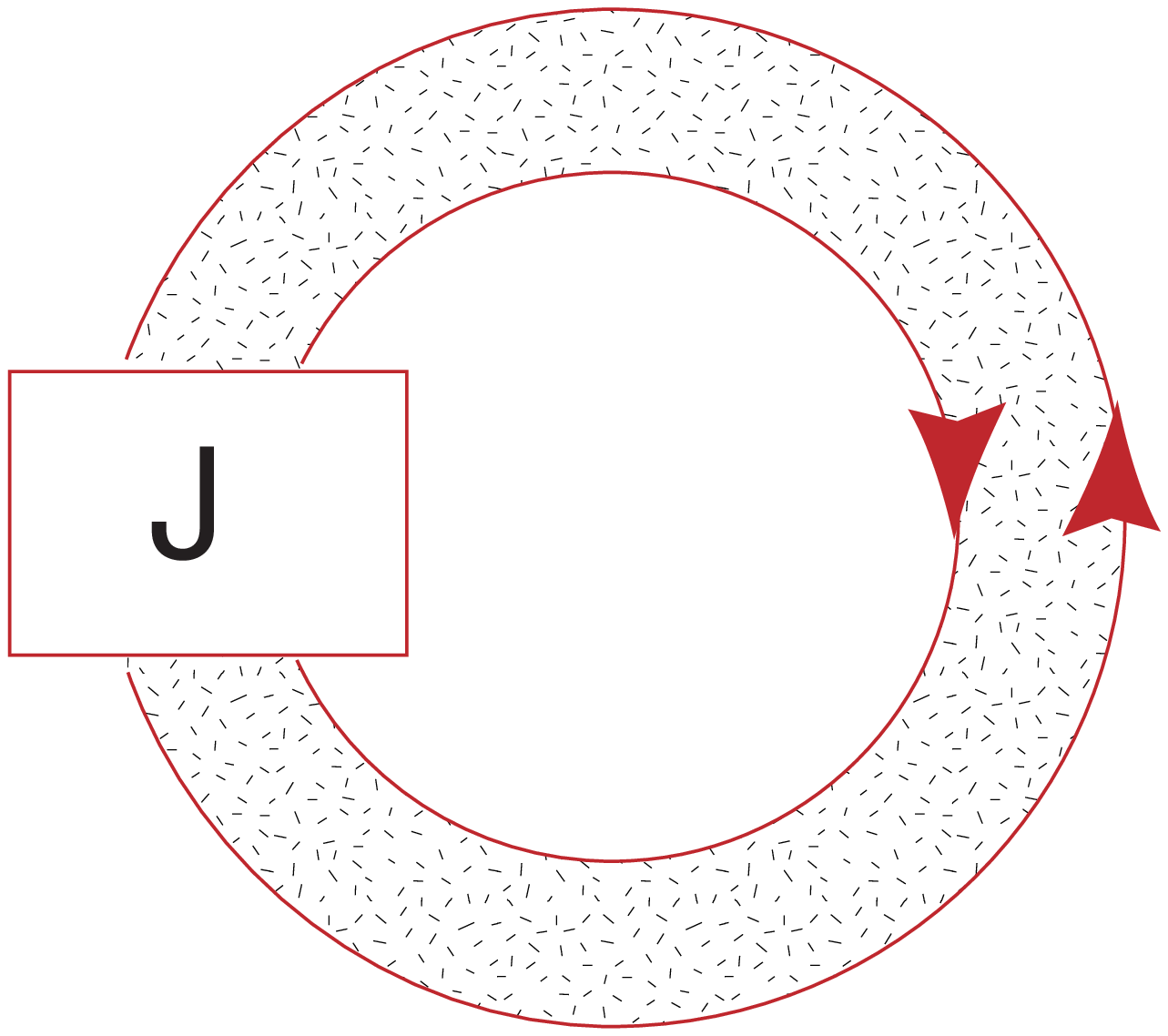}
\hfill
               \end{center}
 \caption{Knotted annuli with oriented sutures, $A(K_{1})$ and $A(K_{2})$} 
\label{fig4}
\end{figure}
{\prop \label{proposition:1}For a given knot $K$ and its meridian $\mu$, we have,
\[ \widehat{HFK}(S^{3}_{N}(K), \mu, g(K)) \cong \left\{
     \begin{array}{lr}
       \widehat{HFK}(S^{3}, K, g(K)) & if N \ne 0\\
       \widehat{HFK}(S^{3},K, g(K)) \oplus \widehat{HFK}(S^{3},K, g(K)) & if N = 0
     \end{array}
   \right . \] }
\begin{proof}
Let $(M, \gamma)$ be the knot complement with two oriented longitudes as the sutures(see Figure \ref{fig4}). We decompose $(M, \gamma)$ along $S(K)$, a Seifert surface of $K$, to obtain $(M^{'}, \gamma^{'})$, 
\[ \xy(0,0)
\xymatrix{(M, \gamma) \ar @{~>} (25,0) ^-{S(K)}="(M, \gamma)" &\text{\space \space \space \space \space \space \space} (M^{'}, \gamma^{'})}
\endxy . \] 
Notice that,
\[ (M^{'}, \gamma^{'}) \cong \left\{
     \begin{array}{lr}
      S^{3}(S(K)) & if N \ne 0\\
      S^{3}(S(K)) \sqcup S^{3} & if N = 0
     \end{array}
   \right . \]
Now by taking $SFH$ of both sides and also from the fact that $SFH(S^{3}(S(K)) \cong \widehat{HFK}(S^{3}, K, g(K))$, the result follows. The last isomorphism is \cite[Theorem 1.5]{Juhasz2008}.
\end{proof}
Proposition \ref{proposition:1}, in fact, gives us a formula for $G_{n}$. Note that $S(K)$ is a nice decomposing surface(see \cite{Juhasz2008} for the definition of a nice decomposing surface) and $(M^{'}, \gamma^{'})$ is taut. Therefore, based on Proposition \ref{prop:1}, $SFH(M^{'}, \gamma^{'})$ is a face of the polytope of $SFH(M, \gamma)$, that is, $SFH(M^{'}, \gamma^{'}) \cong G_{n}$. 
In order to obtain $G_{1}$, we just need to decompose $(M, \gamma)$ along the same surface, $S(K)$, with the opposite orientation. One can also obtain $H_{1}$ and $H_{m}$ in a similar manner by replacing $S(K) = S(K_{1})$ by $S(K_{2})$. Thus, we have computed $G_{1} \otimes H_{1}$, $G_{1} \otimes H_{m}$, $G_{n} \otimes H_{1}$ and $G_{n} \otimes H_{m}$ in the polytope \ref{fig8}. In particular, these four groups are all non-zero since, $\widehat{HFK}(S^{3}, K_{i}, g(K_{i})) \ncong 0$, $i = 1, 2$, provided that $K_{i}$ is nontrivial. \\

Proposition \ref{proposition:1} can be interpreted in the language of \emph{longitude Floer homology} as well. As we mentioned in Remark \ref{rmk:1}, the sutured Floer homology of the Seifert surface complementary manifold corresponding to the knot $K$, is isomorphic to $\widehat{HFL}(K)$. Therefore, when the framing is zero, we obtain the following corollary for $\widehat{HFL}$.
{\cor For a given knot $K$ in $S^{3}$, 
\[ \widehat{HFL}(K, \mathfrak{s}_{top}) \cong  \widehat{HFK}(S^{3},K, g(K)) \oplus \widehat{HFK}(S^{3},K, g(K)) \]
where $\mathfrak{s}_{top}$ is the highest grading in the support of $\widehat{HFL}$,i.e., $\widehat{HFL}(K,  i) \cong 0$ for $ i >  \mathfrak{s}_{top}$.} 
\subsection{Proof of the main theorem}
In this subsection we prove the main theorem of the paper. Having known the results of the subsections \ref{section:3.3} and \ref{section:3.4} , it remains to distinguish $R$ and $R^{'}$. We mention here a natural notion of equivalence for sutured manifolds, where $(M, \gamma)$ and $(M^{'}, \gamma^{'})$ are equivalent if $SFH(M, \gamma)$ and $SFH(M^{'}, \gamma^{'})$ are isomorphic. We still need to clarify the notion of isomorphism between the sutured Floer homology of two sutured manifolds. The following is \cite[Definition 4.1]{Hedden2008}.
{\defn \label{defn:8}Two relatively Spin$^{c}$-graded groups 
\[ SFH(M, \gamma) = \bigoplus_{\mathfrak{s} \in \underline{Spin}^{c}(M, \gamma)} SFH(M, \gamma, \mathfrak{s})  \text{          ,          }SFH(M^{'}, \gamma^{'}) = \bigoplus_{\mathfrak{s} \in \underline{Spin}^{c}(M^{'}, \gamma^{'})} SFH(M^{'}, \gamma^{'}, \mathfrak{s}) \]
are isomorphic if
    \begin{enumerate}
               \item There is an isomorphism, $f^{*} : H^{2}(M^{'}, \partial M^{'}; \mathbb{Z}) \rightarrow  H^{2}(M, \partial M; \mathbb{Z}).$
               \item There is a bijection of sets $u : \underline{Spin}^{c}(M^{'}, \gamma^{'}) \rightarrow  \underline{Spin}^{c}(M, \gamma).$
                \item The following diagram commutes    
\begin{displaymath}
    \xymatrix{
        \underline{Spin}^{c}(M^{'}, \gamma^{'}) \otimes H^{2}(M^{'}, \partial M^{'}; \mathbb{Z})  \ar[r]^{(u, f^{*})} \ar[d] & \underline{Spin}^{c}(M, \gamma) \otimes H^{2}(M, \partial M; \mathbb{Z}) \ar[d] \\
      \underline{Spin}^{c}(M^{'}, \gamma^{'})   \ar[r]_{u}       & \underline{Spin}^{c}(M, \gamma) }
\end{displaymath}
where the vertical arrows are followed by the action of $H^{2}(M^{'}, \gamma^{'})$ on \underline{Spin}$^{c}(M^{'}, \gamma^{'})$ and $H^{2}(M, \gamma)$ on \underline{Spin}$^{c}(M, \gamma).$
\item There are isomorphisms $g_{\mathfrak{s}} : SFH(M^{'}, \gamma^{'}, \mathfrak{s}) \rightarrow  SFH(M, \gamma, u (\mathfrak{s}))$ for every $\mathfrak{s} \in \underline{Spin}^{c}(M^{'}, \gamma^{'})$
\end{enumerate}
}
If $SFH(M, \gamma) \cong SFH(M^{'}, \gamma^{'})$ in the sense of the above definition, $f^{*}$ and $u$ are both obtained by pulling back along a function $f: (M, \gamma) \rightarrow (M^{'}, \gamma^{'})$; an equivalence between $(M, \gamma)$ and $(M^{'}, \gamma^{'})$. In addition, if the surfaces $R$ and $R^{'}$ are weakly equivalent, then $f$ comes from the restriction of $(S^{3}, R) \rightarrow (S^{3}, R^{'})$ to $S^{3}(R)$ that gets mapped to $S^{3}(R^{'})$. Also, $f_{*} : H_{1}(S^{3} \setminus R) \rightarrow H_{1}(S^{3} \setminus R^{'})$ preserves the Seifert form, i.e, $a.b = f_{*}(a).f_{*}(b)$ for every $a,b \in  H_{1}(S^{3} \setminus R)$.
\begin{proof}[Proof of the main theorem] We recall that the Seifert matrices for the two surfaces $R$ and $R^{'}$ are  given by $V_{R} = \left( \begin{array}{cc}
l & 1 \\
0 & 0  \end{array} \right)$  and 
 $V_{R^{'}} = \left( \begin{array}{cc}
l & 0 \\
-1 & 0  \end{array} \right)$ based on the basis elements $c_{i}$'s and $d_{i}$'s, shown in Figure \ref{fig1}($i = 1,2$). Set $M = S^{3}(R)$ and $M^{'} = S^{3}(R^{'})$. Using the notation of Section \ref{section:3}, let us mark one generator in each $G_{i} \otimes H_{j}$ by $x_{ij}$ and let us also denote the group generated by these elements, by $\langle x_{ij} \rangle$. If the two Seifert surfaces were equivalent, then based on Definition \ref{defn:8} we would have $\sigma : SFH(M, \gamma) \rightarrow  SFH(M^{'}, \gamma^{'})$, a bijection from the generators to generators, i.e., it maps every generator of $SFH(M, \gamma)$, say $x_{ij}$, to a generator $\sigma(x_{ij})$ of $ SFH(M^{'}, \gamma^{'})$. We also obtain a map $f: (M, \gamma) \rightarrow (M^{'}, \gamma^{'})$, where $f$ is compatible with taking difference classes, i.e., $\epsilon (a, b) . \epsilon (c, d) = f_{*}\epsilon (a, b) . f_{*}\epsilon (c, d)$ for some isomorphism $f_{*} : H_{1}(S^{3} \setminus R) \rightarrow H_{1}(S^{3} \setminus R^{'})$ where, $a, b, c, d \in H_{1}(S^{3} \setminus R)$.\\
Suppose such a $\sigma$ exists. As we mentioned earlier, in the Seifert forms $V_{R}$ and $V_{R^{'}}$, $l$ is a non-zero integer. Let us assume that $l > 0$. Then,
\[ \epsilon (\sigma (x_{11}), \sigma (x_{nm}))^{2} = (f_{*}\epsilon (x_{11}, x_{nm}))^{2} = (nc_{1} + mc_{2})^{2} = n^{2}l + nm \]
Since $\sigma$ sends generators to generators, for some $\alpha$ and $\beta$ in $\mathbb{Z}$, $\epsilon (\sigma (x_{11}), \sigma (x_{nm})) = \alpha d_{1} + \beta d_{2}$. The fact that $f_{*}$ preserves the Seifert form implies that,
\[ (\alpha d_{1} + \beta d_{2})^{2} =  n^{2}l + nm. \]
On the other hand,
\[  (\alpha d_{1} + \beta d_{2})^{2} = \alpha ^{2}l - \alpha \beta. \]
Notice that $x_{11}$ and  $x_{nm}$ have furthest distance in the polytope of $S^{3}(R)$. This implies that, $n^{2}l + nm$ is the greatest positive number that can possibly
be generated from $\alpha ^{2}l - \alpha \beta$, where $|\alpha | \le n$ and $|\beta | \le m$. Hence, we must must have $|\alpha| = n$ and $|\beta| = m$. In addition, $\alpha$ and $\beta$ have opposite signs. Thus, $\epsilon (\sigma (x_{11}), \sigma (x_{nm})) =\pm( n d_{1} - m d_{2})$. The following illustrates what we just figured out,
\begin{figure}[h]
\begin{center}
\setlength{\unitlength}{0.8cm}
\begin{picture}(8,0)(0,-10)
\put(6.9,-11){\line(1,0){1}}
\put(8.5,-11){\line(1,0){1}}
\put(11.6,-11){\line(1,0){1}}
\put(13.55,-11){\line(1,0){1}}
\put(5.1,-11){$\langle \sigma(x_{ij})\rangle$}
\put(8,-11){$...$}
\put(9.7,-11){$\langle \sigma(x_{k1}) \rangle$}
\put(12.85,-11){$...$}
\put(14.75,-11){$\langle \sigma(x_{11}) \rangle$}
\put(7.2,-10.8){$d_{1}$}
\put(8.9,-10.8){$d_{1}$}
\put(11.95,-10.8){$d_{1}$}
\put(13.9,-10.8){$d_{1}$}
\put(6.6,-11.3){\line(0,-1){1}}
\put(6.6,-13){\line(0,-1){1}}
\put(6.5,-12.6){$.$}
\put(6.5,-12.7){$.$}
\put(6.5,-12.8){$.$}
\put(5.1,-14.6){$\langle \sigma(x_{nm}) \rangle$}
\put(6.1,-11.8){$d_{2}$}
\put(6.1,-13.5){$d_{2}$}
\put(7.1,-14.6){\line(1,0){1}}
\put(8.7,-14.6){\line(1,0){1}}
\put(12,-14.6){\line(1,0){1}}
\put(13.65,-14.6){\line(1,0){1}}
\put(8.2,-14.6){$...$}
\put(9.9,-14.6){$\langle \sigma(x_{km})\rangle$}
\put(13.05,-14.6){$...$}
\put(14.75,-14.6){$\langle \sigma(x_{i^{'}j^{'}} \rangle$}
\put(7.4,-14.4){$d_{1}$}
\put(9.1,-14.4){$d_{1}$}
\put(12.25,-14.4){$d_{1}$}
\put(13.9,-14.4){$d_{1}$}
\put(15.3,-11.3){\line(0,-1){1}}
\put(15.3,-13){\line(0,-1){1}}
\put(15.2,-12.6){$.$}
\put(15.2,-12.7){$.$}
\put(15.2,-12.8){$.$}
\put(9.5,-12.8){$.$}
\put(14.8,-11.8){$d_{2}$}
\put(14.8,-13.5){$d_{2}$}
\put(-4.75,-11){\line(1,0){1}}
\put(-2.8,-11){\line(1,0){1}}
\put(-.4,-11){\line(1,0){1}}
\put(1.55,-11){\line(1,0){1}}
\put(-6,-11){$\langle x_{11}\rangle$}
\put(-3.5,-11){$...$}
\put(-1.6,-11){$\langle x_{k1} \rangle$}
\put(.85,-11){$...$}
\put(2.75,-11){$\langle x_{1n} \rangle$}
\put(-4.40,-10.8){$c_{1}$}
\put(-2.45,-10.8){$c_{1}$}
\put(-.05,-10.8){$c_{1}$}
\put(1.9,-10.8){$c_{1}$}
\put(-5.4,-11.3){\line(0,-1){1}}
\put(-5.4,-13){\line(0,-1){1}}
\put(-5.5,-12.6){$.$}
\put(-5.5,-12.7){$.$}
\put(-5.5,-12.8){$.$}
\put(-6.2,-14.6){$\langle x_{1m} \rangle$}
\put(-5.9,-11.8){$c_{2}$}
\put(-5.9,-13.5){$c_{2}$}
\put(-4.75,-14.6){\line(1,0){1}}
\put(-2.8,-14.6){\line(1,0){1}}
\put(-.2,-14.6){\line(1,0){1}}
\put(1.55,-14.6){\line(1,0){1}}
\put(-3.5,-14.6){$...$}
\put(-1.6,-14.6){$\langle x_{km}\rangle$}
\put(.85,-14.6){$...$}
\put(2.75,-14.6){$\langle x_{nm} \rangle$}
\put(-4.40,-14.4){$c_{1}$}
\put(-2.45,-14.4){$c_{1}$}
\put(-.05,-14.4){$c_{1}$}
\put(1.9,-14.4){$c_{1}$}
\put(3.3,-11.3){\line(0,-1){1}}
\put(3.3,-13){\line(0,-1){1}}
\put(3.2,-12.6){$.$}
\put(3.2,-12.7){$.$}
\put(3.2,-12.8){$.$}
\put(-5.5,-12.8){$.$}
\put(2.8,-11.8){$c_{2}$}
\put(2.8,-13.5){$c_{2}$}
\end{picture}
\end{center}
\end{figure}
\\
\\
\\
\\
\\
\\
\\
\\
\\
 Observe that $\sigma (x_{11})$ and $\sigma (x_{nm})$ are located along the other diagonal compared to $x_{11}$ and $x_{nm}$. Now, in the horizontal direction, take the closest nonzero group $G_{k} \otimes H_{1}$ to $\langle x_{11} \rangle$ such that $k>0$. Such a group exists since, \emph{a priori} $G_{n} \otimes H_{1} \ne 0$. Thus, we can assume that $x_{k1} \ne 0$ for some $k > 0$. Based on our convention, $x_{k1} \in G_{k} \otimes H_{1}$ and $k \le n$. Therefore, we have $\epsilon (x_{11}, x_{k1}) = \pm kc_{1}$. Then,
\[ \epsilon (\sigma (x_{11}), \sigma (x_{k1}))^{2} = (f_{*}\epsilon (x_{11}, x_{k1}))^{2} = k^{2}c_{1}^{2} = k^{2}l, \]
and for some $\alpha$ and $\beta$, $\epsilon (\sigma (x_{11}), \sigma (x_{k1})) = \alpha d_{1} + \beta d_{2}$ which in turn ensures
\[ (\alpha d_{1} + \beta d_{2})^{2} = \alpha ^{2}l - \alpha \beta = k^{2}l, \]
where $k > 0$, and so $\alpha \ne 0$. Since $SFH$ is symmetric and we obtained the polytope from the tensor product formula and also from the way we chose $x_{k1}$, we conclude that $\alpha$ is at least $k$. The way that $\sigma (x_{11})$ is located in the polytope, implies that $\alpha$ and $\beta$ have different signs. Thus, $\alpha \beta \le 0$ and $\alpha^{2} - \alpha \beta \ge \alpha ^{2}$, which implies $\beta = 0$ and $\epsilon (\sigma (x_{11}), \sigma (x_{k1})) =\pm kd_{1}$. We get a contradiction now. On the one hand, 
\[ \epsilon (\sigma (x_{11}), \sigma (x_{nm})) . \epsilon (\sigma (x_{11}), \sigma (x_{k1})) = \pm (nd_{1} - md_{2}) . \pm kd_{1} = \pm (nkl + mk), \]
where $l, m, n, k > 0$. On the other hand, 
\[ \epsilon (x_{11}, x_{nm}) . \epsilon (x_{11}, x_{k1}) = \pm (nc_{1} + mc_{2}) . \pm kc_{1} = \pm nkl . \]
For the case $l < 0$, instead of $x_{11}$, $x_{nm}$ and $x_{k1}$ we take $x_{n1}$, $x_{1m}$ and $x_{(n - k) 1}$ for the smallest positive $k$, where $x_{(n - k) 1} \ne 0$. Then, a contradiction follows similarly. These show that  $SFH(S^{3}(R)) \ncong SFH(S^{3}(R^{'}))$ and so $R \nsimeq R^{'}$.
\end{proof}
\bibliography{References}
\end{document}